\documentclass[11pt]{amsart}

 \usepackage{amsfonts,graphics,amsmath,amsthm,amsfonts,amscd,amssymb,amsmath,latexsym,multicol,
 mathrsfs}
\usepackage{epsfig,url}
\usepackage{flafter}
\usepackage{fancyhdr}
\usepackage{hyperref}
\hypersetup{colorlinks=true, linkcolor=black}
 \usepackage[matrix, arrow]{xy}

\DeclareMathOperator{\Pic}{Pic}

\DeclareMathOperator{\Supp}{Supp}

\DeclareMathOperator{\vol}{vol}

 \numberwithin{equation}{subsection}
 \numberwithin{footnote}{subsection}

 \newtheorem{thm}[subsection]{Theorem}
 \newtheorem{conj}[subsection]{Conjecture}

{
\theoremstyle{upright}

}

 \newcommand{\N}{\mathbb N}
 \newcommand{\PP}{\mathbb P}
 \newcommand{\A}{\mathbb A}
 \newcommand{\Q}{\mathbb Q}
 \newcommand{\R}{\mathbb R}

 \newcommand{\bir}{\dashrightarrow}
 \newcommand{\rddown}[1]{\left\lfloor{#1}\right\rfloor} % round-down

\title{\large G\MakeLowercase{eneralised pairs in birational geometry}}
\author{\large C\MakeLowercase{aucher} B\MakeLowercase{irkar}}
\date{\today}
\begin{document}
\maketitle

\begin{abstract}
In this note we introduce generalised pairs from the perspective of the evolution of the notion of 
space in birational algebraic geometry. We describe some applications of generalised pairs in recent 
years and then mention a few open problems.

\end{abstract}

\tableofcontents

%%%%%%%%%%%%%%%%%%%%%%%%
%%%%%%%%%%%%%%%%%%%%%%%%%%

\bigskip
\bigskip

Unless stated otherwise we work over an algebraically closed field $k$ of characteristic zero. 
Divisors are considered with rational coefficients; for simplicity we avoid real coefficients.
Although we will make some historical remarks but this text is not meant to tell the history of birational geometry. 
This note is an expanded version of a talk given at the conference ``Geometry at large" in 
Fuerteventura, Spain, in December 2018. 

Thanks to Yifei Chen for many helpful comments.

\section{Varieties} 

Geometry is the study of ``shapes", that is, ``spaces". The notion of space varies in different types of geometries 
and it also changes historically while the subject evolves over time. The spaces of interest in 
classical algebraic geometry are algebraic varieties $X$ defined in affine or projective spaces 
by vanishing of polynomial equations over an algebraically closed field $k$. It is not 
easy to classify such spaces even in dimension one due to the huge number of possibilities. 
It is then natural to restrict to smaller classes of varieties such as smooth ones which 
leads one to birational geometry.

Birational geometry aims to classify algebraic varieties up to birational transformations. This involves 
finding special representatives in each birational class and then classifying such representatives, say 
by constructing their moduli spaces. To start with, taking compactification 
and then normalisation one can assume $X$ to be projective and normal. However, being normal 
is not a very special property so it is not easy to work in such generality. 
We can resolve the singularities 
and assume $X$ to be smooth. This enables us to use many tools available for smooth varieties, e.g. 
the Riemann-Roch theorem in low dimension has a particularly simple expression. 

In dimension one each birational class contains a unique smooth projective model 
so smoothness is a satisfactory condition to work with. One then can focus on parametrising 
curves of fixed genus as did Riemann in the 19th century.

In dimension two one can still stick with the category of smooth varieties. 
Indeed running a classical minimal model program on $X$ by successively blowing down 
$-1$-curves, smoothness is preserved. The program ends with either $\PP^2$ or a $\PP^1$-bundle 
over a curve or a smooth projective variety $Y$ with $K_Y$ being nef. Here nefness means that 
the intersection number $K_Y\cdot C$ is non-negative for every curve $C\subset Y$. However, we are 
still led to consider some kind of singularities. Indeed if $K_Y$ is big, then some positive multiple  of
it is base point free hence defines a birational contraction $Y\to Z$ where $Z$ is a projective surface with 
Du Val singularities. The two-dimensional birational geometry was developed by the 
Italian school of algebraic geometry in late 19th century and early 20th century, in particular, by 
Castelnuovo and Enriques.

It is clear that the primary spaces to study in dimension one and two are smooth 
projective varieties. Many of the deep statements in classical algebraic geometry are stated for 
such varieties. However, later when people studied moduli spaces, say of curves of fixed genus, 
in greater detail they were led to consider more general 
spaces such as curves with nodal singularities which may not be irreducible [\ref{Deligne-Mumford}]. 
Constructing moduli spaces for surfaces involves considering 
semi-lc surfaces which are analogues of curves with nodal curves (c.f. [\ref{kollar-moduli}][\ref{Alexeev}]).

Until the end of the 1960's there were some progress in higher dimensional birational geometry, 
e.g. work of Fano, resolution of singularities of 3-folds by Zariski, and 
resolution of singularities in arbitrary dimension by Hironaka. 
It was only in the 1970's that the subject really took off. For example, Iskovskikh and 
Manin [\ref{Isk-Manin}] proved that smooth quartic 3-folds are not rational using earlier work of Fano.  
On the other hand, Iitaka proposed a program for birational classification of varieties 
(even open varieties) according to their Kodaira dimension (c.f. [\ref{Iitaka-1}][\ref{Iitaka-2}]) which as 
we will see led to the development of pairs. And the techniques developed in Mori's solution to Hartshorne 
conjecture [\ref{Mori-Hartshorne-conj}] in late 1970's  paved the way for taking some of the first steps toward a 
minimal model program in higher dimension [\ref{Mori-non-nef}]. 
In yet another direction Reid's study of singularities [\ref{Reid}] 
helped to identify the right classes of singularities for this program to work. 
Moreover, Kawamata and Shokurov [\ref{Kawamata-fg}][\ref{Shokurov-nonvanishing}] made extensive use of pairs 
and vanishing theorems [\ref{Kawamata-van}] to establish foundational results 
such as the base point free theorem guaranteeing existence of relevant contractions generalising 
contraction of $-1$-curves on surfaces. With Mori's proof of existence of flips [\ref{Mori-flip}],
(and the more general version by Shokurov [\ref{Shokurov-log-flips}]), and 
Miyaoka and Kawamata's proof of abundance [\ref{Miyaoka-1}][\ref{Miyaoka-2}][\ref{KMM-abundance}] 
(and the more general version 
by Keel, Matsuki, M$\rm ^c$Kernan [\ref{KMM-abundance}]), and further input from Koll\'ar and many others  
birational geometry for 3-folds was well-developed by the early 1990's.
See [\ref{KMM}][\ref{kollar-mori}] for introductions to the subject.

As it is evident from the above discussion,  
it did not take long before people understood that the category of smooth varieties is not large 
enough for birational geometry in dimensions $\ge 3$. One needs to consider singularities albeit 
of mild kinds. There are several type of singularities such as Kawamata log terminal and log canonical 
which behave well (these are defined in Section 3). These singularities not only allow proofs to go through but their study 
reveals fascinating facts which are often reflections of global phenomena.   
Thus the primary spaces in higher dimensional birational geometry were no longer smooth ones but 
those with mild singularities.

On the other hand, one can study the birational geometry of varieties $X$ 
defined over arbitrary fields $k$. In characteristic zero many 
statements can be reduced to similar statements after passing to the algebraic closure of $k$, e.g. 
running a minimal model program. But the finer classification is not that simple, for example, 
there can be many non-isomorphic smooth projective varieties over $k$ which all become 
isomorphic after passing to the algebraic closure. See for example [\ref{Kollar-3-folds/R}] for birational geometry 
of varieties of dimension $3$ over $\R$.   
In positive characteristic the story is more complicated because of non-separable  
field extensions.

\section{Schemes}
With Grothendieck's revolution of algebraic geometry attention shifted from studying varieties 
to schemes in the 1960's, to some extent. Birational geometry of arithmetic surfaces, that is, 
regular schemes of dimension two was worked out [\ref{Shafarevich}][\ref{Lichtenbaum}]. 
This is important for arithmetic geometry. e.g. for 
constructing proper regular models of curves defined over number fields and eventually 
to construct N\'eron models of elliptic curves [\ref{Silverman}]. 

On the other hand, one may attempt to study birational geometry of schemes defined over a field; 
we mean allowing non-reduced structures. 
There does not seem to be much research in this direction except some cases in dimension one. 
For example, Mumford's canonical curves 
appear while studying minimal smooth projective surfaces $Y$ of Kodaira dimension one [\ref{Mumford-surfaces-char-p}]. 
An effective divisor $X$ on $Y$ such that $K_Y\cdot X=0$ and $X\cdot C=0$ for 
every component $C$ of $X$, is a curve of canonical type. By adjunction 
$K_X\sim 0$ so we may consider $X$ as a Calabi-Yau scheme of dimension one.    
Studying such $X$, especially when $X$ is connected and the g.c.d of its coefficients 
is $1$, gives non-trivial information about $X$.

From the 1970's focus was on the birational geometry of good old smooth 
varieties but one important point of view of Grothendieck remained which emphasized on working in the 
relative setting, that is, studying varieties with a morphism to another variety. 
This has had an important influence on the development of birational geometry. 
For example a Kawamata log terminal singularity 
can be considered as a local analogue of a Fano variety.

\section{Pairs}

In the 1970's Iitaka initiated a program to study the birational classification of open varieties (that is, non-compact 
varieties) according to their Kodaira dimension [\ref{Iitaka-1}][\ref{Iitaka-2}]. This involved compactifying 
the variety by adding a boundary divisor. This approach then evolved into the theory of pairs. 
Pairs are roughly speaking algebraic varieties together with an extra structure given by a certain type of 
divisor (see below for precise definition). 
They appear naturally even if one is only interested in studying projective varieties. We give a 
series of motivating examples before defining pairs rigorously. 

Much of the machinery of birational geometry of the last four decades has been developed 
for pairs. The main concept of space is then that of a pair.

\subsection{Open varieties}  
Suppose $U$ is a smooth variety. Then the geometry of $U$ is often best understood after taking a 
compactification. Suppose $X$ is a smooth projective variety containing $U$ as an open subset. Such a 
compactification is not unique except in dimension one. We can choose $X$ so that $B:=X\setminus U$ 
is a divisor with simple normal crossing singularities. Now we consider $(X,B)$ as a pair. 
The geometry of $(X,B)$ reveals so much about the geometry of $U$. 

\subsection{Adjunction}
Suppose $X$ is a smooth variety and $B$ is a smooth prime divisor on $X$. 
Then the well-known adjunction formula says that $K_B=(K_X+B)|_B$. Generalisations of the 
adjunction formula play a central role in birational geometry. Very often one derives 
non-trivial statements about the geometry of $(X,B)$ from that of $B$ allowing 
proofs by induction on dimension. Thus one is led to study pairs such as $(X,B)$.

\subsection{Canonical bundle formula}
Suppose $X$ is a smooth projective variety and $f\colon X\to Z$ is a contraction, that is, a 
projective morphism with connected fibres. Suppose in addition that $K_X\sim f^*L$ for some 
$\Q$-divisor $L$. Then the canonical bundle formula says that we can can write 
$$
K_X\sim f^*(K_Z+B+M)
$$
where the discriminant divisor $B$ is a $\Q$-divisor uniquely determined and the moduli divisor 
$M$ is a $\Q$-divisor determined up to $\Q$-linear equivalence. The classical example 
is that of Kodaira's canonical bundle formula in which $X$ is a surface and $f$ is an elliptic fibration. 
The higher dimensional version in a more general context is discussed below.
The canonical bundle formula allows one to investigate the geometry of $X$ from that of $(X,B+M)$. 
It is possible to choose $M$ so that we can consider $(X,B+M)$ as a pair. 
We will see later that it is actually more appropriate to consider it as a generalised pair.   

\subsection{Quotient varieties}
Suppose that $X$ is a smooth variety and $G$ is a finite group acting on $X$. 
Let $Y$ be the quotient of $X$ by $G$ and $\pi\colon X\to Y$ be the quotient map. 
Then using Hurwitz formula one can write $K_X=\pi^*(K_Y+B_Y)$ for some divisor $B_Y$ 
whose coefficients are of the form $1-\frac{1}{n}$ for certain natural numbers $n$. 
It is then natural to study $Y$ and $B_Y$ together rather than just $Y$, that is, 
to consider the pair $(Y,B_Y)$.

\subsection{Definition of pairs}
We now define pairs and their singularities rigorously.
A \emph{pair} $(X,B)$ consists of a normal quasi-projective variety $X$ and a $\Q$-divisor 
$B\ge 0$ such that $K_X+B$ is $\Q$-Cartier. We call $B$ a \emph{boundary} divisor if its coefficients 
are in $[0,1]$. 
Let $\phi\colon W\to X$ be a log resolution of $(X,B)$. We can write 
$$
K_W+B_W=\phi^*(K_X+B)
$$ 
where $B_W$ is uniquely determined. Here we choose $K_W$ so that $\phi_*K_W=K_X$. 
We say $(X,B)$ is \emph{log canonical} (resp. \emph{Kawamata log terminal})(resp. \emph{$\epsilon$-log canonical}) if 
every coefficient of $B_W$ is $\le 1$ (resp. $<1$)(resp. $\le 1-\epsilon$). 
If $D$ is a component of $B_W$ with coefficient $1$, then its image on $X$ is called a 
\emph{log canonical centre}.

\subsection{Examples}
(1)
The simplest kind of pairs are the log smooth ones. 
A \emph{log smooth} pair is a pair $(X,B)$ where $X$ is smooth and $\Supp B$ has simple 
normal crossing singularities. Such pairs are log canonical. If the coefficients of $B$ are 
less than $1$, then the pair is Kawamata log terminal.

(2) 
Let $X=\PP^2$. When $B$ is a nodal curve, then $(X,B)$ is log canonical. But 
when $B$ is a cuspidal curve, then $(X,B)$ is not log canonical.

(3)
Let $X$ be the cone over a smooth rational curve. Then $(X,0)$ is Kawamata log terminal. 
In contrast if $Y$ is the cone over an elliptic curve, then $(Y,0)$ is log canonical but 
not Kawamata log terminal.

(4)
Let $X$ be a toric variety and $B$ be the sum of the torus invariant divisors. 
Then $(X,B)$ is log canonical.

(5)
Let $X$ be the variety in $\A^4$ defined by the equation $xy-zw=0$. Then 
$(X,0)$ is Kawamata log terminal.

\subsection{Using pairs}
We illustrate the power of pairs by an example which frequently comes up in inductive statements.
Suppose $(X,B)$ is a projective Kawamata log terminal pair and $S$ is a normal prime divisor on $X$. 
Now suppose $L$ is a Cartier divisor such that 
$$
A:=L-(K_X+B+S)
$$ 
is ample. Consider the exact sequence 
$$
H^0(X,L)\to H^0(S,L|_S)\to H^1(X,L-S).
$$
By assumption, 
$$
L-S=K_X+B+A. 
$$
Since $A$ is ample and $(X,B)$ is Kawamata log terminal, $H^1(L-S)=0$ by the Kawamata-Viehweg vanishing theorem 
(c.f. [\ref{Kawamata-van}]). Therefore, every section of $L|_S$ can be lifted to a section of $L$. 
This is very useful for example in situations where we want to show the linear system $|L|$ is 
non-empty or that it is base point free.

\subsection{Progress in the last two decades}
There has been huge progress in birational geometry of pairs in arbitrary dimension, 
in the last two decades which builds 
upon the machinery developed in the preceding decades. This is due to work of 
many people. To name few specific examples, existence of flips 
[\ref{Shokurov-pl}][\ref{HM-flip}][\ref{BCHM}], existence of minimal models for varieties of log general 
type [\ref{BCHM}], ACC for log canonical thresholds and boundedness results for 
varieties of log general type [\ref{HMX1}][\ref{HMX2}] (ACC on smooth varieties was proved 
in [\ref{dFEM}]), boundedness of complements and Fano 
varieties [\ref{B-compl}][\ref{B-BAB}], existence of moduli spaces for stable pairs [\ref{kollar-moduli}], 
existence of moduli spaces for polarised Calabi-Yau and Fano pairs [\ref{B-moduli-pol-var}]  
have all been established.

\subsection{Semi-log canonical pairs}
A semi-log canonical pair is roughly a pair in which the underlying set may not be irreducible 
or may have components with self-intersections but the intersections should be nice similar to those on nodal curves. 
Semi-log canonical pairs are important for constructing compact moduli spaces [\ref{kollar-moduli}] 
as they appear as limits of usual pairs in families.

%%%%%%%%%%%%%%%%%%%%%%%%%%%%%

\section{Generalised pairs}

In recent years a new concept of space has evolved, that is, the generalised pairs. 
A generalised pair is roughly a pair together with a birational model polarised with some divisor 
having some positivity property. They were first defined and studied in [\ref{BZh}] in their general form.  
However, some special cases were already investigated in [\ref{BH}] and simpler forms also appeared  
implicitly in the earlier work [\ref{B-WZD}]. 
Generalised pairs have found applications in various contexts 
which we will discuss in the next section. 
We begin with considering some motivating examples of such pairs and then give their precise definition afterwards.

\subsection{Polarised varieties}
Consider a projective variety $X$ and an ample divisor $M$ on it. We say $X$ is polarised 
by $M$. For example, 
$M$ can be a very ample divisor determining an embedding of $X$ into some projective space. 
Polarised varieties play a central role in moduli theory as one often needs some kind of 
positivity in order to achieve ``stability" hence be able to construct moduli spaces. 
For example, one may consider varieties $X$ polarised by $K_X$. 

For applications it is important to allow more general polarisations, that is, when $M$ 
is not necessarily ample but only nef. Indeed, [\ref{BH}] studies the birational geometry of 
varieties and pairs polarised by nef divisors which paved the way for the development of 
generalised pairs. Here $M$ is only a divisor class and not 
necessarily effective so we cannot consider $(X,M)$ as a usual pair. 

Moduli spaces of polarised Calabi-Yau and Fano pairs were constructed recently in 
[\ref{B-moduli-pol-var}].

\subsection{Generalised polarised varieties}
Suppose $X$ is a projective variety and $X\bir X'$ is a birational map to another projective variety.
Suppose $X'$ is polarised by a nef divisor $M'$. We can consider $X$ as a generalised polarised 
variety, that is, a variety with some polarised birational model. Sometimes it is important 
to understand how $X$ and $X'$ are related. For example, when $X'$ is a  minimal model of 
$X$, then the canonical divisor $K_{X'}$ polarises $X'$ so we can consider $X$ as a generalised 
polarised variety. This can be used to prove results about termination of flips (see the next section).

\subsection{Canonical bundle formula and subadjunction}
Suppose $(X,B)$ is a projective pair with Kawamata log terminal singularities and $f\colon X\to Z$ is a contraction. 
Suppose in addition that $K_X+B\sim_\Q f^*L$ for some 
$\Q$-divisor $L$. Then the canonical bundle formula says that we can can write 
$$
K_X+B\sim_\Q f^*(K_Z+B_Z+M_Z)
$$
where the discriminant divisor $B_Z$ is a $\Q$-divisor uniquely determined and the moduli divisor 
$M_Z$ is a $\Q$-divisor determined up to $\Q$-linear equivalence [\ref{kaw-subadjuntion}]. There is a 
resolution $Z'\to Z$ and a nef divisor $M_{Z'}$ on $Z'$ whose pushdown to $Z$ is $M_Z$ [\ref{ambro-lc-trivial}]. 
We can then consider $(Z,B_Z+M_Z)$ as a generalised pair which happens to be generalised Kawamata log terminal 
according to the definitions below. 

Now assume $(Y,\Delta)$ is a projective log canonical pair and $V$ is a minimal log canonical centre. 
Here minimality is among log canonical centres with respect to inclusion. 
Assume $(Y,\Theta)$ is Kawamata log terminal for some $\Theta$. Assume $V$ is normal (this 
actually holds automatically). Then we can write 
$$
(K_Y+\Delta)|_V\sim_\Q K_V+C+N
$$
where $C$ is a divisor with coefficients in $[0,1]$ and $N$ is a divisor class which is the pushdown of some 
nef divisor. What happens is that from the setup we can find a Kawamata log terminal pair 
$(S,\Gamma)$ and a contraction $g\colon S\to V$ such that 
$$
K_S+\Gamma\sim_\Q g^*(K_Y+\Delta)|_V.
$$
Thus we can use the previous paragraph to decompose $(K_Y+\Delta)|_V$ as $K_V+C+N$. 
In particular, $(V,C+N)$ is a generalised pair. This construction is called \emph{subadjunction} 
[\ref{kaw-subadjuntion}] - in practice though people often try to perturb $C+N$ to get a boundary divisor.

\subsection{Definition of generalised pairs}
A projective \emph{generalised pair} consists of 
\begin{itemize}
\item a normal projective variety $X$, 

\item a $\Q$-divisor $B\ge 0$ on $X$, and 

\item a birational contraction $\phi\colon X'\to X$ and a nef $\Q$-divisor $M'$ on $X'$ 
\end{itemize}
such that $K_{X}+B+M$ is $\Q$-Cartier where $M:= \phi_*M'$. 

Actually we specify $X',M'$ only up to birational transformations, that is, 
if we replace $X'$ with a resolution and replace $M'$ with its pullback, then the pair would be the same.
We can similarly define generalised pairs in the relative setting, that is, when 
we are given a projective morphism $X\to Z$ (but $X$ may not be projective) 
and then we only assume $M'$ to be nef over $Z$.

Now we define generalised singularities for a generalised pair $(X,B+M)$.
Replacing $X'$ we can assume $\phi$ is a log resolution of $(X,B)$. We can write 
$$
K_{X'}+B'+M'=\phi^*(K_{X}+B+M)
$$
for some uniquely determined $B'$. 
We say $(X,B+M)$ is 
\emph{generalised log canonical} (resp. \emph{generalised Kawamata log terminal})
(resp. \emph{generalised $\epsilon$-log canonical}) if each coefficient of $B'$ is $\le 1$ 
(resp. $<1$)(resp. $\le 1-\epsilon$).

\subsection{Examples}
We present a series of examples of generalised pairs.\\
(1)
The most obvious way to construct a generalised pair is to take a projective pair $(X,B)$ and a
nef $\Q$-divisor $M$ to get $(X,B+M)$.   
Then $(X,B+M)$ is generalized log canonical (resp. generalized Kawamata log terminal) iff $(X,B)$ is
log canonical (resp. Kawamata log terminal). In this example $M'=M$ does not
contribute to the singularities even if its coefficients are large. In contrast, the larger
the coefficients of $B$, the worse the singularities.

(2)
In general, $M'$ does contribute to singularities. For example, assume $X=\PP^2$ and that $\phi$ is
the blowup of a point $x$. Let $E'$ be the exceptional divisor, $L$ a line passing through
$x$ and $L'$ the birational transform of $L$.

If $B=0$ and $M'=2L'$, then
we can calculate $B'=E'$ hence $(X,B+M)$ is generalized log canonical but not generalized Kawamata log terminal.
However, if $B=L$ and $M'=2L'$, then
$(X,B+M)$ is not generalized log canonical because in this case $B'=L'+2E'$.

(3) 
Suppose that $X$ is a smooth projective variety such that $-K_X$ is nef. 
Letting $M:=-K_X$ we can consider $(X,M)$ as a generalised Calabi-Yau pair as $K_X+M=0$. 
This is useful for tackling certain problems, see \ref{ss-bnd-lcy} below.

(4) 
Let $X$ be a smooth projective variety and $D$ be a divisor with $h^0(X,D)\neq 0$. 
We can find a resolution of singularities $\phi\colon X'\to X$ such that the 
movable part $M'$ of the linear system $|\phi^*D|$ is base point free, hence in particular nef. 
Writing $C'$ for the fixed part of $|\phi^*D|$ we get $M'+C'\sim \phi^*D$. 
Let $M,C$ be the pushdowns of $M',C'$ respectively. Then $(X,C+M)$ is a generalised pair 
with nef part $M'$ which remembers the movable and fixed parts of $|\phi^*D|$. Its 
singularities reflect the singularities of $(X,L)$ where $L$ is a general member of 
$|D|$, e.g. $(X,C+M)$ is generalised log canonical iff $(X,L)$ is log canonical. 
Moreover, for any non-negative rational number $t$, the singularities of the generalised pair $(X,tC+tM)$ with nef part $tM'$ 
reflect the singularities of the pair $(X,tL)$.

More generally suppose $W\subseteq |D|$ is a linear system. Again pick a resolution 
$X'\to X$ so that the 
movable part $M'$ of the linear system $\phi^*W$ is base point free, hence in particular nef. 
Writing $C'$ for the fixed part of $\phi^*W$ we get $M'+C'\sim \phi^*D$.
Let $M,C$ be the pushdowns of $M',C'$ respectively. Then $(X,C+M)$ is a generalised pair with nef part $M'$ 
which remembers the movable and fixed parts of $W$. Similarly, for a non-negative rational number $t$, 
$(X,tC+tM)$ is a generalised pair with nef part $tM'$ whose behaviour reflects that of 
$(X,tL)$ where $L$ is a general member of $W$. Generalised pairs of this kind have 
appeared in the context of birational rigidity and non-rationality of Fano varieties and 
the Sarkisov program going back to work of Iskovskikh-Manin which in turn is based on ideas 
going back as far as Fano and Neother (see [\ref{Cheltsov}] for a survey on these topics, especially section 0.2).

For simplicity we assumed $X$ to be smooth but similar constructions apply for example if $X$ 
has klt singularities and $D$ is $\Q$-Cartier. 

(5)
Let $X$ be a smooth projective variety such that $K_X\sim_\Q 0$ (that is, $X$ is a Calabi-Yau variety). 
Let $M$ be a nef divisor on $X$. A difficult conjecture predicts that $M$ is numerically equivalent 
to a semi-ample divisor. Viewing $(X,M)$ as a generalised pair, the question is whether $K_X+M$ 
is numerically semi-ample. See the last section for more general statements.

(6)
Let $X$ be a normal $\Q$-factorial variety and let $V=\sum b_i V_i$ be a formal linear combination of 
closed subvarieties where $V_i\neq X$ and $b_i$ are non-negative rational numbers. One can consider $(X,V)$ as some kind of 
pair generalising the traditional notion of a pair [\ref{EMY}] (when all the $V_i$ are prime divisors, then $(X,V)$ is a 
pair in the usual sense). Take a log resolution 
$\phi\colon X'\to X$ such that the scheme-theoretic inverse images $V_i':=\phi^{-1}V_i$ are all Cartier 
divisors. Here we are assuming that the exceptional locus of $\phi$ union the support of all the $V_i'$ 
is a divisor with simple normal crossing singularities (hence the name log resolution). 
Write $K_{X'}+E'=\phi^*K_X$ and then define the log discrepancy of a prime divisor 
$D'$ on $X'$ with respect to $(X,V)$ to be $1-\mu_{D'}(E'+\sum b_iV_i')$. We say $(X,V)$ is log canonical 
(resp. Kawamata log terminal) if 
all the log discrepancies are $\ge 0$ (resp. $>0$). 

The singularities of $(X,V)$ can be interpreted in terms of generalised pairs. 
Indeed, for each $i$, the morphism $X'\to X$ factors through the blowup $Y_i\to X$ of $X$ 
along $V_i$. The scheme-theoretic inverse image of $V_i$ on $Y_i$ is a Cartier divisor which is anti-ample 
over $X$, and $V_i'$ is the pullback of this anti-ample divisor. Thus $-V_i'$ is nef over $X$, for each $i$. 
Now let $M':=-\sum b_iV_i'$ which is nef over $X$, and let $B':=E'+\sum b_iV_i'$. Consider the generalised pair 
$(X,B+M)$ with nef part $M'$ relatively over $X$ where $B=\phi_*B'$ and $M=\phi_*M'$. Then we have 
$$
K_{X'}+B'+M'=K_{X'}+E'=\phi^*K_X=\phi^*(K_X+B+M)
$$ 
as $B+M=0$.
Thus the log discrepancies of $(X,V)$ are the same as the generalised log discrepancies of 
$(X,B+M)$. 

Note that even though $B+M=0$ but $(X,B+M)$ as a generalised pair is not the same as the 
usual pair $(X,0)$. Also note that even if $X$ is projective, the generalised pair $(X,B+M)$ 
only makes sense relatively over $X$ since $M'$ obviously may not be nef globally.
This suggests that $(X,B+M)$ is useful for studying local properties of $(X,V)$.

\subsection{Geometry of generalised pairs}
In the next section we present some of the applications of generalised pairs in recent years. 
We want to emphasize that in addition to such applications studying the geometry of 
generalised pairs is interesting on its own. Many questions for varieties and usual pairs can be 
asked in the context of generalised pairs which leads to some deep problems. For example, 
can we always run a minimal model program for a generalised pair and get a minimal model 
or a Mori fibre space in the end? See the last section for some more problems. 

Various geometric aspects of generalised pairs have been studied in recent years, 
for example, see [\ref{BH}][\ref{BZh}][\ref{Moraga-term}][\ref{Hacon-Moraga}] for the minimal model program 
and termination, [\ref{BZh}] for birational boundedness of linear systems and ACC for generalised 
lc thresholds, [\ref{B-compl}] for boundedness of complements, [\ref{LP-1}][\ref{LP-2}] 
for abundance, [\ref{Filipazzi}] for boundedness of generalised pairs of general type, 
[\ref{Filipazzi-adjunction}] for the canonical bundle formula, 
[\ref{liu-sarkisov}] for the Sarkisov program, [\ref{B-moduli-pol-var}] for birational boundedness 
of linear systems and boundedness of polarised varieties, [\ref{liu-lc-thresholds}] for 
accumulation points of generalised lc thresholds,
[\ref{filipazzi-svaldi}] for invariance of plurigenera and boundedness, 
[\ref{FMX}][\ref{Filipazzi-Moraga}][\ref{Chen}] for more on boundedness of complements. 

%%%%%%%%%%%%%%%%%%%%%%

\section{Some applications of generalised pairs}

In this section we discuss some applications of generalised pairs in recent years.

\subsection{Effective Iitaka fibrations and pluricanonical systems of generalised pairs}
Let $W$ be a smooth projective variety of Kodaira dimension $\kappa(W) \ge 0$. The 
Kodaira dimension $\kappa(W)$ is the largest number $\kappa$ among $\{-\infty, 0, 1,\dots,\dim X\}$ 
such that 
$$
\limsup_{m\in \N} \frac{h^0(X,mK_W)}{m^\kappa}>0.
$$ 
Then by a well-known construction of Iitaka, there is a birational
morphism $V\to W$ from a smooth projective variety $V$, and a contraction $V\to X$ onto
a projective variety $X$ such that a (very) general fibre $F$ of $V\to X$
is smooth with Kodaira dimension zero, and $\dim X$ is equal to the Kodaira dimension
$\kappa(W)$. The map $W\bir X$ is referred to
as an \emph{Iitaka fibration} of $W$, which is unique up to birational equivalence.
For any sufficiently divisible natural number $m$, the pluricanonical
system $|mK_W|$ defines an Iitaka fibration. The following hard conjecture predicts that 
we can choose $m$ uniformly depending only on the dimension. 

\begin{conj}[Effective Iitaka fibration, {\rm c.f. [\ref{HM}]}]\label{conj-eff-iitaka-fib}
Let $W$ be a smooth projective variety of dimension $d$ and Kodaira dimension $\kappa(W)\ge 0$.
Then there is a natural number $m_d$ depending only on $d$ such that the pluricanonical system
$|mK_W|$ defines an Iitaka fibration for any natural number $m$ divisible by $m_d$.
\end{conj}

In [\ref{BZh}] the conjecture is reduced to bounding certain invariants of the very general fibres of 
the Iitaka fibration. The point is that bounding such invariants, perhaps after replacing $W,X$ with birationally, 
the Iitaka fibration induces a canonical bundle type formula giving $\Q$-divisors $B\ge 0$ and $M$  
where the coefficients of $B$ are in a DCC set and $M$ is nef with bounded Cartier index. Thus it would 
be enough to prove the following theorem regarding generalised pairs. Here by DCC set we mean 
the set does not contain any infinite strictly decreasing sequence of numbers.

\begin{thm}\label{t-bir-bnd-M}
Let $\Lambda$ be a DCC set of nonnegative real numbers,
and $d,r$ be natural numbers. Then there is a natural number $m(\Lambda, d,r)$
depending only on $\Lambda, d,r$ such that if:

\begin{itemize}
\item[(i)]
$(X,B)$ is a projective log canonical pair of dimension $d$,
\item[(ii)]
the coefficients of $B$ are in $\Lambda$,
\item[(iii)]
$rM$ is a nef Cartier divisor, and
\item[(iv)]
$K_X+B+M$ is big,
\end{itemize}
then the linear system $|{m(K_X+B+M)}|$ defines a birational map if $m\in \N$ is divisible by $m(\Lambda, d,r)$.
\end{thm}

For usual pairs, that is when $M=0$, the theorem was previously known  [\ref{HMX2}, Theorem 1.3].
However, the general case proved in [\ref{BZh}] uses very subtle properties of the theory of generalised 
pairs. Note that for a $\Q$-divisor $D$, by $|D|$ and $H^0(X,D)$ we mean  $|\rddown{D}|$ and
$H^0(X,\rddown{D})$.\\

\subsection{Boundedness of complements and of Fano varieties} 
Let $(X,B)$ be a projective pair where $B$ is a boundary. 
Let $T=\rddown{B}$ and $\Delta=B-T$. 
An \emph{$n$-complement} of $K_{X}+B$ is of the form 
$K_{X}+{B}^+$ such that 
\begin{itemize}
\item $(X,{B}^+)$ is log canonical, 

\item $n(K_{X}+{B}^+)\sim 0$, and 

\item $n{B}^+\ge nT+\rddown{(n+1)\Delta}$.
\end{itemize}
From the definition one sees that 
$$
-nK_{X}-nT-\rddown{(n+1)\Delta}\sim n{B}^+-nT-\rddown{(n+1)\Delta}\ge 0
$$
so existence of an $n$-complement for $K_X+B$ implies that the linear system 
$$
|-nK_{X}-nT-\rddown{(n+1)\Delta}|
$$
is non-empty. In particular, this means that we should be looking at varieties $X$ with 
$K_X$ ``non-positive", e.g. Fano varieties. Complements were defined by Shokurov [\ref{Shokurov-log-flips}] 
in the context of construction of flips. 
The following was conjectured by him [\ref{shokurov-surf-comp}] and 
proved in [\ref{B-compl}].

\begin{thm}\label{t-bnd-compl-usual}
Let $d$ be a natural number and $\mathfrak{R}\subset [0,1]$ be a finite set of rational numbers.
Then there exists a natural number $n$ 
depending only on $d$ and $\mathfrak{R}$ satisfying the following.  
Assume $(X,B)$ is a projective pair such that 
\begin{itemize}

\item $(X,B)$ is log canonical of dimension $d$,

\item the coefficients of $B$ are in $\Phi(\mathfrak{R})$, 

\item $X$ is of Fano type, and 

\item $-(K_{X}+B)$ is nef.\\
\end{itemize}
Then there is an $n$-complement $K_{X}+{B}^+$ of $K_{X}+{B}$ 
such that ${B}^+\ge B$. Moreover, the complement is also an $mn$-complement for any $m\in \N$. 
\end{thm}

Here $\Phi(\mathfrak{R})$ stands for the set  
$$
\left\{1-\frac{r}{m} \mid r\in \mathfrak{R},~ m\in \N\right\}.
$$
A special case of the theorem is when $K_X+B\sim_\Q 0$ along a fibration $f\colon X\to T$. 
This is where generalised pairs come into the picture. 
Applying the {{canonical bundle formula}} we can write 
$$
K_X+B\sim_\Q f^*(K_T+B_T+M_T)
$$
where $B_T$ is the {discriminant divisor} and $M_T$ is the {moduli divisor}. It turns out that 
the coefficients of $B_T$ are in $\Phi(\mathfrak{S})$ 
for some fixed finite set $\mathfrak{S}$ of rational numbers, and that $pM_T$ is integral for some bounded number $p\in \N$. 
Now we want to find a complement for $K_T+B_T+M_T$ and pull it back to $X$. 
As mentioned elsewhere in the text $(T,B_T+M_T)$ is not a pair but it is a 
generalised pair. Thus we actually need to construct complements 
 in the more general setting of generalised pairs which can be defined similar to the case of usual pairs. 
 Once we have a bounded complement for 
$K_T+B_T+M_T$ we pull it back to get a bounded complement for $K_X+B$. 

The theory of complements is applied in [\ref{B-BAB}] to prove the following statement which was 
known as the BAB conjecture. Thus generalised pairs play an important (indirect) role in the proof of this theorem.

\begin{thm}
Let $d$ be a natural number and $\epsilon$ be a positive real number. Then the projective 
varieties $X$ such that  

$\bullet$ $(X,B)$ is $\epsilon$-log canonical of dimension $d$ for some boundary $B$, and 

$\bullet$  $-(K_X+B)$ is nef and big,\\
form a bounded family. 
\end{thm}

This theorem in turn has been applied to various problems in recent years.

\subsection{Termination of flips and existence of minimal models} 
In [\ref{B-WZD}] existence of minimal models is linked with existence of weak forms of Zariski decompositions. 
A given divisor $D$ on a projective variety has a weak Zariski decomposition if its pullback 
to some resolution of $X$ can be written as $P+N$ where $P$ is nef and $N$ is effective. When 
$(X,B)$ is a pair, we would be interested in weak Zariski decompositions of $K_X+B$. 
On the other hand, termination of flips is linked with log canonical thresholds in [\ref{B-term-acc}]. 
Extending these to the case of generalised pairs, termination of flips 
for generalised pairs with weak Zariski decompositions is derived from termination in lower dimension  
for generalised pairs, in [\ref{Hacon-Moraga}][\ref{Han-Li}]; it is also shown that 
existence of weak Zariski decompositions for pseudo-effective 
generalised pairs is equivalent to existence of minimal models for such pairs.  
In particular, termination of flips is established for pseudo-effective pairs of dimension four 
which at the moment does not follow from any other technique (this was first established in [\ref{Moraga-term}] 
for usual pairs).

\subsection{Boundedness of certain rationally connected varieties}\label{ss-bnd-lcy}
M${\rm ^c}$Kernan and Prokhorov [\ref{McP}] conjectured a more general form of BAB. 

\begin{conj}
Let $d$ be a natural number and $\epsilon$ be a positive real number. Consider 
projective varieties $X$ such that 
\begin{itemize}
\item $(X,B)$ is $\epsilon$-log canonical of dimension $d$ for some boundary $B$, 

\item $-(K_X+B)$ is nef, and 

\item $X$ is rationally connected. 
\end{itemize}
Then the set of such $X$ forms a bounded family. 
\end{conj}

The rational connectedness assumption cannot be removed: indeed it is well-known that 
K3 surfaces do not form a bounded family; they 
satisfy the assumptions of the theorem with $d=2$ and $\epsilon=1$ and $B=0$ except that 
they are not rationally connected. 

The conjecture fits nicely into the framework of generalised pairs and related 
conjectures. Indeed letting $M:=-(K_X+B)$ we get a generalised $\epsilon$-log canonical  pair $(X,B+M)$ 
with $K_X+B+M= 0$, hence a generalised Calabi-Yau pair. The advantage of this point of view 
 is that running a minimal model program preserves the Calabi-Yau condition, hence one can get 
information by passing to Mori fibre spaces when $B+M\not\equiv 0$ (which is always the case 
on some birational model). Using this and the machinery developed in [\ref{B-lcyf}], a slightly weaker form of 
the conjecture is verified in dimension three in [\ref{BDS}] where one replaces boundedness with boundedness 
up to isomorphism in codimension one.

\subsection{Varieties fibred over abelian varieties}
In [\ref{B-JC}], generalised pairs, more precisely, polarised pairs, are used to investigate 
pairs that are relatively of general type over a variety of maximal Albanese dimension.   
Suppose that $(X, B)$ is a projective Kawamata log terminal pair and 
$f \colon X \to Z$ is a surjective morphism where $Z$ is a normal projective variety
with maximal Albanese dimension, e.g. an abelian variety. It is shown that if $K_X + B$ is big 
over $Z$, then $(X, B)$ has a good log minimal model. Moreover, if $F$ is a general fibre of $f$, then
$$
\kappa(K_X + B) \ge \kappa(K_F + B_F ) + \kappa(Z) = \dim F + \kappa(Z)
$$
where $K_F + B_F = (K_X + B)|_F$.

%%%%%%%%%%%%%%%%%%%%%%
\section{Problems}
In this section we discuss some open problems regarding generalised pairs.

\subsection{Existence of contractions and flips} 
Let $(X,B+M)$ be a projective generalised log canonical pair. Assume $R$ is a $K_X+B+M$-negative extremal ray. 
When $(X,B+M)$ is generalised Kawamata log terminal, existence of the contraction associated to $R$ follows easily 
from the similar result for usual pairs. This is because we can easily find an ample $\Q$-divisor $A$ 
and a Kawamata log terminal pair $(X,\Delta)$ such that 
$$
K_X+\Delta\sim_\Q K_X+B+M+A
$$
and such that $R$ is $(K_X+\Delta)$-negative. In particular, if $R$ defines a flipping contraction, 
then its flip exists [\ref{BCHM}]. 

Now assume $(X,B+M)$ is not generalised Kawamata log terminal. If $(X,C)$ is Kawamata log terminal 
for some $C$, then we can take an 
average and use the previous paragraph. More precisely, taking $t>0$ to be a small rational number,
$$
(X,tC+(1-t)B+(1-t)M)
$$ 
is generalised Kawamata log terminal and 
$$
K_X+tC+(1-t)B+(1-t)M
$$ 
intersects $R$ negatively. Thus in this case the contraction of $R$ exists and in the flipping case 
its flip exists. If there is no $C$ as above, e.g. when $X$ itself has some non-Kawamata log terminal singularities, 
then the situation is more complicated and as far as we know 
existence of contractions and flips in full generality is not proved yet. This is important for running 
minimal model program for generalised pairs.

\subsection{Generalised minimal model program}
Suppose that $(X,B+M)$ is a projective generalised log canonical pair. 
Assuming that existence of contractions and flips are established for such pairs, 
we can run the minimal model program on $K_X+B+M$ 
which, if terminates, produces a generalised Mori fibre space or a generalised minimal model. 
Termination of the program does not seem to follow from termination for usual pairs (although 
we can inductively treat the case when $K_X+B+M$ is pseudo-effective, as in [\ref{Hacon-Moraga}]). 
For polarised varieties it was proved in [\ref{BH}] that termination of \emph{some} 
choice of minimal model program can be guaranteed if termination holds for usual pairs which is only known 
up to dimension three. In low dimension the picture is clearer. 
Generalised termination holds trivially in dimension two. It is also known in dimension 
three [\ref{Moraga-term}] and in dimension four in the pseudo-effective case [\ref{Hacon-Moraga}].

\subsection{Generalised abundance}
Although generalised pairs behave like usual pairs in many ways but there are 
some crucial differences. Assume $(X,B+M)$ is a projective generalised Kawamata log terminal pair with 
$K_X+B+M$ nef. In the case of usual pairs, that is, when $M=0$, the abundance conjecture 
says that $K_X+B$ is semi-ample which means that $m(K_X+B)$ is base point free for some natural number $m$. 
In general when $M\neq 0$, we cannot expect $K_X+B+M$ to be semi-ample as was already pointed out 
in [\ref{BH}]. Indeed this already fails in dimension 
one when $X$ is an elliptic curve, $B=0$, and $M$ is a numerically trivial but non-torsion 
divisor. In dimension one at least numerical abundance holds, that is, $K_X+B+M\equiv D$ for some 
semi-ample divisor $D$. In dimension two even numerical abundance fails. The obstructions 
to numerical abundance seem to arise only when $K_X+B$ is not pseudo-effective (in which case 
$X$ is uniruled so rational curves play a role); it is conjectured that 
if $K_X+B$ is pseudo-effective and if $M$ is nef, then numerical abundance holds [\ref{LP-1}][\ref{LP-2}] (these papers 
prove some results in this direction).

\subsection{Boundedness of generalised pairs with fixed volume} 
Fix natural numbers $d,p,v$ and a DCC set $\Phi\subset [0,1]$ of rational numbers. 
Consider projective generalised log canonical pairs $(X,B+M)$ of dimension $d$ 
such that the coefficients of $B$ are in $\Phi$, 
the divisor $pM'$ is Cartier where $M'$ is the nef part of the pair $(X,B+M)$, 
and $K_X+B+M$ is ample with volume $(K_X+B+M)^d=v$. Then it is expected that such $X$ form a bounded family. 
That is, one expects to find a natural number $l$ depending only on $d,r,v,\Phi$ 
such that $l(K_X+B+M)$ is very ample hence defining an embedding of $X$ into 
some fixed projective space $\PP^n$ such the image of $X$ under this embedding has 
bounded degree in $\PP^n$. This statement is verified in dimension two in [\ref{Filipazzi}] but otherwise 
seems to be open.

\subsection{Classification of generalised pairs}
One may ask to classify generalised pairs in given dimension. 
The extra information in generalised pairs as opposed to usual pairs makes the classification 
theory richer and more complex. For example let's look at the case of dimension one. Let $(X,B+M)$ be a 
projective generalised log canonical pair of dimension one. To classify such generalised pairs 
we need to fix some invariants. We can fix the degree of $K_X+B+M$ and 
also assume that $lB$, $lM$ are both integral divisors for some fixed natural number $l$. 
In particular, the number of components of $B$ is bounded. In fact, $(X,B)$ varies 
in some bounded family of pairs. On the other hand, $lM$ is a nef Cartier divisor whose degree 
takes only finitely many possibilities but $(X,\Supp (B+M))$ need not belong to a bounded family.  
Indeed even when $X=\PP^1$, $B=0$, and $\deg M>0$, then $\Supp M$ can have any number of 
components. One idea is to consider 
$(X,B+M)$ up to some kind of equivalence, e.g. up to isomorphism of $(X,B)$ and up to linear 
equivalence of $lM$.  For example if we fix $X$ to be an elliptic 
curve and $B=0$ and $M\equiv 0$ and then consider $lM$ up to linear equivalence, 
then the classes of such $(X,M)$ are parametrised by the elements 
of $\Pic^0(X)$ (by sending $(X,M)$ to $lM$). Varying both $X$ and $M$ is parametrised by a much larger space.
It is not hard to imagine that going to higher dimension the classification 
problem gets quite complicated.

Another idea is to consider only the case $M\ge 0$. Indeed even in higher dimension this 
works in some interesting situations. For example, consider projective semi-log canonical Calabi-Yau pairs 
$(X,B)$ of fixed dimension $d$ and ample Weil divisors $M\ge 0$ such that $cB$ is a Weil divisor 
for some fixed rational number $c>0$, 
the volume $\vol(M)=v$ is fixed, and $(X,B+tM)$ is log canonical for some $t>0$ (here $t$ is not fixed). 
Then it is shown in [\ref{B-moduli-pol-var}, Theorem 1.10] that there is a projective coarse moduli space for such 
$(X,B),M$ with the fixed data $d,c,v$. Note that $M$ being ample, $(X,B+M)$ is a generalised pair.

%%%%%%%%%%%%%%%%%%%%%%%%%%%%%%%%%%%%%

\vspace{2cm}
%%%%%%%%%%%%%%%%%%%%%

\textsc{DPMMS, Centre for Mathematical Sciences} \endgraf
\textsc{University of Cambridge,} \endgraf
\textsc{Wilberforce Road, Cambridge CB3 0WB, UK} \endgraf
\email{c.birkar@dpmms.cam.ac.uk\\}

\end{document}